\documentclass[reqno,11pt,a4paper]{amsart}

\usepackage{amsthm,amsmath,amssymb}
\usepackage{mathrsfs,amsfonts,mathtools}
\usepackage[colorlinks=true,linkcolor=blue,citecolor=blue,urlcolor=blue]{hyperref}

\newtheorem{theorem}{Theorem}[section]
\newtheorem{proposition}[theorem]{Proposition}

\newtheorem{corollary}[theorem]{Corollary}
\newtheorem{remark}[theorem]{Remark}
\newtheorem{lemma}[theorem]{Lemma}

\numberwithin{equation}{section}
\allowdisplaybreaks
\newcommand{\ddbar}{\sqrt{-1}\,\partial\bar\partial}
\newcommand{\tr}{\operatorname{tr}}
\newcommand{\osc}{\operatorname{osc}}
\newcommand{\cH}{\mathcal H}
\newcommand{\cC}{\mathcal C}
\newcommand{\cQ}{\mathcal Q}
\newcommand{\dosc}{d_{\mathrm{osc}}}

\begin{document}

\title[Stability of a logarithmic inverse-trace flow]{Comparison principles and dynamical stability for a logarithmic inverse-trace flow on compact Hermitian manifolds}

\author[L. Zhang]{Liangdi Zhang}
\address[Liangdi Zhang]{Mathematics Science Research Center\\Chongqing University of Technology\\Chongqing 400054, P. R. China}
\email{ldzhang91@163.com}

\begin{abstract}
We study a logarithmic inverse-trace flow on compact connected Hermitian
manifolds. We first establish a parabolic comparison principle and the
resulting oscillation nonexpansiveness for admissible potentials modulo
constants. We then introduce common-drift sub- and supersolutions.
Their combination yields a uniform drift-corrected zero-order estimate for
arbitrary admissible initial data, while the lower barrier alone implies the
strict cone condition required in the second-order estimate. Combining these
observations with Sun's Hermitian second-order estimate, parabolic
Evans--Krylov regularity, Schauder estimates, and a uniform Harnack inequality,
we obtain long-time existence and exponential smooth convergence. In
particular, whenever the associated elliptic equation admits a smooth
solution, its class is the unique globally attracting fixed point of the
induced semiflow, and any two trajectories synchronize exponentially modulo
constants. We also derive a Hermitian cone--supersolution convergence
criterion and an affine-renormalization criterion that uses Sun's convergence
theorem to construct one elliptic solution and then propagates global
dynamical stability to every admissible initial potential.


\end{abstract}

\maketitle

\section{Introduction}

Let $M$ be a compact connected complex manifold of complex dimension
$n\geq2$, and let $\omega$ and $\chi$ be Hermitian metrics on $M$.  Set
\[
\mathcal H_\chi
=
\left\{
\varphi\in C^\infty(M,\mathbb R):
\chi_\varphi:=\chi+\sqrt{-1}\partial\bar\partial\varphi>0
\right\}.
\]
For $0<T\leq\infty$, by an admissible solution on $[0,T)$ we mean a function
\[
\varphi\in C^\infty(M\times[0,T),\mathbb R)
\]
such that $\varphi(\cdot,t)\in\mathcal H_\chi$ for every $t\in[0,T)$.

We use $[\chi]$ to denote the affine $\partial\bar\partial$-class
\[
[\chi]
=
\left\{
\chi+\sqrt{-1}\partial\bar\partial f:
f\in C^\infty(M,\mathbb R)
\right\},
\]
without assuming that $\chi$ is closed.

For $F\in C^\infty(M,\mathbb R)$, define
\begin{equation}\label{eq:Q}
\mathcal Q_F(\varphi)
=
\log\frac{\chi_\varphi^n}
{\omega\wedge\chi_\varphi^{n-1}}+F
=
\log n-\log\operatorname{tr}_{\chi_\varphi}\omega+F.
\end{equation}
We consider
\begin{equation}\label{eq:flow}
\begin{cases}
\displaystyle
\frac{\partial\varphi}{\partial t}
=
\mathcal Q_F(\varphi),
\\[4pt]
\varphi(\cdot,0)=\varphi_0\in\mathcal H_\chi.
\end{cases}
\end{equation}
The associated elliptic problem is to find
$(\varphi_\infty,\lambda)\in\mathcal H_\chi\times\mathbb R$ such that
\begin{equation}\label{eq:elliptic}
\mathcal Q_F(\varphi_\infty)=\lambda,
\end{equation}
or equivalently
\begin{equation}\label{eq:elliptic-wedge}
\omega\wedge\chi_{\varphi_\infty}^{n-1}
=
e^{F-\lambda}\chi_{\varphi_\infty}^n.
\end{equation}

For constant $F$, the stationary equation associated with \eqref{eq:flow} is the
inverse-trace form of Donaldson's $J$-equation, whereas the logarithmic flow
\eqref{eq:flow} is different from the classical $J$-flow introduced by
Donaldson and Chen \cite{don99,xchen00,xchen04}.  The $J$-flow and related
fully nonlinear inverse and quotient flows have been studied extensively; see
\cite{we04,we2,sw08,flsw14,ls15,cs17,dk19,hk18,gchen2021,flm11,fl12,sun19,tdt22}
and the references therein.  The complex Monge--Amp\`ere equation on compact Hermitian manifolds was
established by Tosatti--Weinkove \cite{tw10}, while its parabolic counterpart
and related regularizing phenomena were studied by Gill \cite{mg11} and
T\^o \cite{to18}.  The parabolic Monge--Amp\`ere theory was subsequently
extended to the almost Hermitian setting by Chu \cite{chu18}. For inverse and quotient equations on compact Hermitian manifolds, Li \cite{yli14}
established elliptic estimates for Donaldson's equation.  Sun \cite{sun16}
subsequently proved elliptic existence, up to a multiplicative constant, under
a strict cone condition together with an additional pointwise inequality.  He
also proved long-time existence for a broad class of parabolic complex
Monge--Amp\`ere type equations; on general Hermitian backgrounds, convergence
requires an additional pointwise condition on the initial data
\cite{sun15}.  Equation \eqref{eq:flow} is the inverse-trace case
$\alpha=1$ of this flow.  Building on the elliptic $\mathcal C$-subsolution
framework \cite{sz18}, Phong and T\^o \cite{pt21} developed a general theory of
fully nonlinear parabolic equations on compact Hermitian manifolds.  Related
Hermitian Hessian quotient flows were studied in \cite{zc22}.

The purpose of this paper is to isolate, for the inverse-trace operator, a
direct comparison-based mechanism which removes the initial-speed restriction
for arbitrary admissible initial data once one elliptic solution, or more
generally a pair of common-drift barriers, is available.  Thus the emphasis is
on global dynamical stability and oscillation nonexpansiveness, rather than on a
claim that the cone condition alone implies convergence.  This distinction is
important in the non-K\"ahler Hermitian setting: the fully Hermitian solvability
criterion stated below retains the additional global supersolution used in the
known elliptic and parabolic existence results.  The key analytical mechanism
is the order-preserving structure of the inverse-trace operator.  The operator in \eqref{eq:Q} is monotone
with respect to the Hermitian matrix order; hence stationary sub- and
supersolutions, after addition of the same linear function of time, trap every
trajectory.  This gives the zero-order estimate directly, without a
space--time Moser iteration.  The lower barrier also implies the strict cone
condition.  Sun's second-order estimate, complex parabolic Evans--Krylov theory,
Schauder estimates, and a uniform parabolic Harnack inequality then yield
long-time existence and exponential convergence.

The relation with \cite{sun15,pt21} is as follows.  We emphasize first that
long-time existence under the cone condition is already contained in Sun's
Theorem~1.1 \cite{sun15}; it is used here as part of the existing analytic
background and is not claimed as a new contribution.  For convergence on a
general Hermitian manifold, Sun's Theorem~1.2 additionally assumes the
pointwise initial inequality
\[
 \frac{\chi_{\varphi_0}^{n}}
 {\omega\wedge\chi_{\varphi_0}^{n-1}}
 \leq e^{-F},
\]
after absorbing the initial potential into the background form.  Equivalently,
this is the sign condition $\cQ_F(\varphi_0)\leq0$.  The bounded case of the
Phong--T\^o theory likewise requires an additional monotonicity or zero-order
input \cite{pt21}.  Our barrier theorem does not assert that these hypotheses
can simply be discarded.  Instead, for the inverse-trace operator, a pair of
sub- and supersolutions with the same drift produces the missing uniform
zero-order estimate by comparison, and hence gives exponential convergence
from every admissible initial potential once such barriers are available.  The
affine criterion at the end of the paper uses Sun's sign condition only once,
to construct one elliptic solution; the comparison theorem then propagates
stability to all other initial data.  We do not claim that the strict cone
condition by itself is sufficient on an arbitrary Hermitian manifold.

There is a second point concerning the zero-order estimate.  The elliptic
argument of Tosatti--Weinkove \cite{twasian} integrates over the spatial
manifold $M$.  Treating time merely as a parameter does not, by itself, control
the quantity $v(x,t)-\inf_{M\times[0,t]}v$ appearing in Sun's parabolic
second-order estimate, since the latter contains information from all earlier
times.  We therefore do not use an elliptic Moser iteration to obtain a
parabolic $C^0$ bound.  Instead, the common-drift comparison estimate gives a
uniform two-sided bound for the shifted potential on the whole space--time
cylinder.  Consequently the space--time infimum in Sun's estimate is controlled
directly.  This point is made quantitatively in Lemma~\ref{lem:sun-quantity}
below.

A related comparison-and-barrier mechanism has recently been developed by
Xiong, Yang, and Yau for prescribed Hermitian--Yang--Mills flows
\cite{xyy26a}.  Their sequel obtains uniform control by a matrix-valued
Uhlenbeck--Yau blow-up procedure and a destabilizing coherent subsheaf
\cite{xyy26b}.  The order-preserving aspect has a precise scalar counterpart in
the present problem.  By contrast, a scalar potential has no nontrivial
spectral projection, so the coherent-subsheaf construction of \cite{xyy26b}
does not directly apply here.

For $H\in C^\infty(M,\mathbb R)$, define the cone condition
\begin{equation}\label{eq:cone}
 \cC_1(H)=\left\{[\chi]:\begin{array}{l}
 \text{there exists }\widehat\chi\in[\chi],\ \widehat\chi>0,\text{ such that}\\[2pt]
 n\widehat\chi^{n-1}>(n-1)e^{-H}\widehat\chi^{n-2}\wedge\omega
 \end{array}\right\}.
\end{equation}
This is the case $\alpha=1$ of the cone condition in \cite{sun15}.
Our first theorem is the basic common-drift barrier principle.

\begin{theorem}[Common-drift barrier principle]\label{thm:barrier}
Assume that there exist $\underline\varphi,\overline\varphi\in\cH_\chi$ and
$\lambda\in\mathbb R$ such that
\begin{equation}\label{eq:barriers}
 \cQ_F(\underline\varphi)\geq\lambda,
 \qquad
 \cQ_F(\overline\varphi)\leq\lambda
 \quad\text{on }M.
\end{equation}
Then, for every $\varphi_0\in\mathcal H_\chi$, equation
\eqref{eq:flow} admits a unique smooth solution on
$M\times[0,\infty)$.  There exists
\[
\varphi_\infty=\varphi_\infty(\varphi_0)\in\mathcal H_\chi
\]
satisfying
\begin{equation}\label{eq:limit-barrier}
 \cQ_F(\varphi_\infty)=\lambda
\end{equation}
such that, for every integer $k\geq0$, there are constants
$C_k,\eta_k>0$ for which
\begin{equation}\label{eq:exp-barrier}
 \|\varphi(t)-\lambda t-\varphi_\infty\|_{C^k(M)}
 \leq C_ke^{-\eta_kt}.
\end{equation}
The class $[\varphi_\infty]\in\mathcal H_\chi/\mathbb R$ is independent
of the initial potential and is the unique solution class of the elliptic
equation.
\end{theorem}

Taking one elliptic solution as both barriers gives the principal dynamical
consequence.

\begin{corollary}[Global dynamical stability]\label{cor:stability}
Suppose that \eqref{eq:elliptic} has a solution $(\varphi_*,\lambda)$.  Then, from
every $\varphi_0\in\cH_\chi$, the solution of \eqref{eq:flow} exists for all
time and $\varphi(t)-\lambda t$ converges exponentially in $C^\infty$ to a
solution of \eqref{eq:elliptic}.  In particular, if $\varphi$ and $\psi$ are
two solutions with arbitrary admissible initial values, then for some constant
$c_\infty$ and every $k\geq0$,
\begin{equation}\label{eq:synchronization-intro}
 \|\varphi(t)-\psi(t)-c_\infty\|_{C^k(M)}
 \leq C_ke^{-\eta_kt}.
\end{equation}
\end{corollary}

The barrier formulation is equivalent to elliptic solvability and global
convergence, but it is useful because it identifies the two independent
analytic inputs: the lower barrier supplies the cone condition, while the pair
of barriers supplies the drift-corrected zero-order estimate.

\begin{corollary}[Equivalent elliptic and dynamical conditions]\label{cor:equivalence}
Fix $\lambda\in\mathbb R$.  The following statements are equivalent:
\begin{enumerate}
 \item Equation \eqref{eq:elliptic} has a solution with constant $\lambda$.
 \item There exist $\underline\varphi,\overline\varphi\in\cH_\chi$ satisfying
 \eqref{eq:barriers}.
 \item For every $\varphi_0\in\cH_\chi$, the solution of \eqref{eq:flow}
exists globally and $\varphi(t)-\lambda t$ converges in $C^\infty$ to an
element of $\cH_\chi$.
\end{enumerate}
Whenever these conditions hold, the convergence in (3) is exponential.
\end{corollary}

Finally, the following sufficient condition produces the elliptic solution
needed in Corollary \ref{cor:stability}.  It records precisely what follows from
an affine time renormalization and Sun's convergence theorem.

\begin{theorem}[Affine-renormalization criterion]\label{thm:affine}
Assume that there exist $\varphi_*\in\cH_\chi$ and $\kappa\in\mathbb R$ such that
\begin{equation}\label{eq:affine-speed}
 \kappa\geq\max_M\cQ_F(\varphi_*)
\end{equation}
and
\begin{equation}\label{eq:affine-cone}
 [\chi]\in\cC_1(F-\kappa).
\end{equation}
Then there is a unique constant $\lambda\in\mathbb R$ for which
\eqref{eq:elliptic} is solvable.  Consequently, for every
$\varphi_0\in\cH_\chi$, the solution of \eqref{eq:flow} exists globally and
$\varphi(t)-\lambda t$ converges exponentially in $C^\infty$ to a solution of
\eqref{eq:elliptic}.
\end{theorem}

The case $\kappa=0$ gives a convenient formulation which involves only the
original, unshifted equation.

\begin{corollary}[Hermitian cone--supersolution criterion]\label{cor:hermitian}
Assume that
\begin{equation}\label{eq:hermitian-cone}
 [\chi]\in\cC_1(F)
\end{equation}
and that there exists $\varphi_*\in\cH_\chi$ satisfying
\begin{equation}\label{eq:hermitian-super}
 \chi_{\varphi_*}^n
 \leq e^{-F}\,\omega\wedge\chi_{\varphi_*}^{n-1},
\end{equation}
or equivalently $\cQ_F(\varphi_*)\leq0$ on $M$.  Then there is a unique
constant $\lambda\in\mathbb R$ for which \eqref{eq:elliptic} is solvable, and
the corresponding solution is unique modulo constants.  For every
$\varphi_0\in\cH_\chi$, the solution of \eqref{eq:flow} exists for all time and,
for every integer $k\geq0$, there are constants $C_k,\eta_k>0$ such that
\begin{equation}\label{eq:hermitian-exp}
 \|\varphi(t)-\lambda t-\varphi_\infty\|_{C^k(M)}
 \leq C_ke^{-\eta_kt}.
\end{equation}
Here $\varphi_\infty$ denotes the drift-corrected limit determined by the
initial potential.
\end{corollary}

\begin{remark}\label{rem:cone-alone}
Condition \eqref{eq:hermitian-cone} alone is not asserted to imply solvability
in the general non-K\"ahler Hermitian setting.  The additional supersolution
\eqref{eq:hermitian-super} is precisely the extra global hypothesis appearing
in Sun's Hermitian existence theorem \cite[Theorem~1.2]{sun16}.  The new
consequence here is that, after one elliptic solution has been obtained, the
comparison principle yields exponential convergence from every admissible
initial potential.
\end{remark}

A pointwise sufficient condition for \eqref{eq:affine-cone} is
\begin{equation}\label{eq:pointwise-cone}
 \chi_{\varphi_*}-\frac{n-1}{n}e^{\kappa-F}\omega>0.
\end{equation}
The initial-speed condition is used only to obtain one elliptic solution from
Sun's theorem.  Once that solution is known, Corollary \ref{cor:stability}
gives convergence for all admissible initial data.  Corollary
\ref{cor:hermitian} is the unshifted fully Hermitian form of this criterion.

The paper is organized as follows.  Section \ref{sec:prelim} records the
algebraic structure, concavity, short-time theory, and the evolution of the
speed.  Section \ref{sec:comparison} proves comparison and nonexpansiveness on
$\cH_\chi/\mathbb R$.  Section \ref{sec:cone} develops common-drift barriers and
the cone condition.  Sections \ref{sec:estimates} and \ref{sec:convergence}
establish the uniform estimates and exponential convergence.  The barrier
results, the fully Hermitian cone--supersolution consequence, and the
affine-renormalization criterion are completed in Sections
\ref{sec:barrier-proof} and \ref{sec:affine-proof}.

\section{Algebraic structure and short-time theory}\label{sec:prelim}

We first record the algebraic identities and evolution equations which will be
used throughout the paper.  For every positive Hermitian $(1,1)$-form $A$, one has
\begin{equation}\label{eq:wedge-trace}
 \omega\wedge A^{n-1}=\frac{1}{n}(\tr_A\omega)A^n.
\end{equation}
Consequently, if
\[
 G(A):=\log\frac{A^n}{\omega\wedge A^{n-1}},
\]
then
\begin{equation}\label{eq:G-trace}
 G(A)=\log n-\log\tr_A\omega.
\end{equation}

In local holomorphic coordinates, write $A=(A_{i\bar j})$, its inverse by
$(A^{i\bar j})$, and $$\omega=\sqrt{-1}\,g_{i\bar j}
 dz^i\wedge d\bar z^j.$$  If $B=(B_{i\bar j})$ is a real $(1,1)$-form, direct
differentiation gives
\begin{equation}\label{eq:DG-expanded}
 DG_A(B)=a_A^{i\bar j}B_{i\bar j},
 \qquad
 a_A^{i\bar j}:=
 \frac{A^{i\bar q}A^{p\bar j}g_{p\bar q}}{\tr_A\omega}.
\end{equation}
The Hermitian matrix $(a_A^{i\bar j})$ is positive definite.  Thus $G$ is
strictly increasing with respect to the Hermitian matrix order.

We shall also use the concavity of $G$.  Let $A_s=A+sB$,  for $s$ in an
interval on which $A_s>0$, and set
\[
 S=\tr(A^{-1}g),\qquad
 P=\tr(A^{-1}BA^{-1}g),\qquad
 R=\tr(A^{-1}BA^{-1}BA^{-1}g).
\]
Since
\[
\frac{d}{ds}A_s^{-1}\bigg|_{s=0}
=-A^{-1}BA^{-1},
\]
we have
\[
S'(0)=-P,\qquad S''(0)=2R.
\]
As $G(A_s)=\log n-\log S(s)$,
\begin{equation}\label{eq:G-second}
 \frac{d^2}{ds^2}G(A_s)\bigg|_{s=0}
 =\frac{P^2}{S^2}-\frac{2R}{S}.
\end{equation}

To verify the sign, put
\[
C=A^{-1/2}BA^{-1/2},
\qquad
D=A^{-1/2}gA^{-1/2}.
\]
Then
\[
S=\tr D,\qquad P=\tr(CD),\qquad R=\tr(C^2D).
\]
Choose a unitary basis in which
$D=\operatorname{diag}(d_1,\ldots,d_n)$, with $d_i>0$.  Since $C$ is
Hermitian,
\[
S=\sum_i d_i,\qquad
P=\sum_i d_i C_{i\bar i},\qquad
R=\sum_{i,j}d_i|C_{i\bar j}|^2.
\]
Therefore the Cauchy--Schwarz inequality yields
\[
P^2
\leq
\left(\sum_i d_i\right)
\left(\sum_i d_i|C_{i\bar i}|^2\right)
\leq SR.
\]
Substituting this into \eqref{eq:G-second}, we obtain
\begin{equation}\label{eq:G-concave}
 D^2G_A(B,B)\leq-\frac{R}{S}\leq0.
\end{equation}
In particular, $G$ is concave on the cone of positive Hermitian forms.

For an admissible solution $\varphi$ of \eqref{eq:flow}, the linearized operator
is
\begin{equation}\label{eq:Lt}
 L_t v=a^{i\bar j}(t)\partial_i\partial_{\bar j}v,
 \qquad
 a^{i\bar j}(t)=
 \frac{\chi_\varphi^{i\bar q}\chi_\varphi^{p\bar j}g_{p\bar q}}
 {\tr_{\chi_\varphi}\omega}.
\end{equation}
For each admissible $\varphi$, the coefficient matrix
$(a^{i\bar j})$ is positive definite. Standard fully nonlinear
parabolic theory gives the following local result.

\begin{proposition}[Short-time existence and uniqueness]\label{prop:short-time}
For every $\varphi_0\in\cH_\chi$, there is a unique smooth admissible solution
of \eqref{eq:flow} on a maximal interval $[0,T_{\max})$, where
$0<T_{\max}\leq\infty$.
\end{proposition}

\begin{proof}
Since $\chi_{\varphi_0}>0$ and $M$ is compact, the eigenvalues of
$\chi_{\varphi_0}$ with respect to $\omega$ are bounded above and away from
zero.  Formula \eqref{eq:DG-expanded} therefore shows that the linearization
at $\varphi_0$ is uniformly elliptic.  By continuity, the same is true in a
sufficiently small $C^2$-neighborhood of $\varphi_0$.  Since $G$ is smooth on
the cone of positive Hermitian forms, the standard local theory for uniformly
parabolic fully nonlinear scalar equations gives a smooth admissible solution
for a short time.  Uniqueness follows from the same local theory, or from the
comparison principle proved in Proposition \ref{prop:comparison}.  The usual
continuation procedure then gives a maximal existence interval
$[0,T_{\max})$.
\end{proof}

Let
\begin{equation}\label{eq:u-def}
 u:=\partial_t\varphi.
\end{equation}
Since $F$ is independent of time, differentiating \eqref{eq:flow} in time and using \eqref{eq:DG-expanded} gives
\begin{equation}\label{eq:u-evolution}
 \partial_tu=L_tu.
\end{equation}

Indeed
\[
 \partial_tu
 =DG_{\chi_\varphi}(\sqrt{-1}\,\partial\bar\partial u)
 =a^{i\bar j}u_{i\bar j}.
\]
For every $T<T_{\max}$, the equation is uniformly parabolic on
$M\times[0,T]$.  Hence the parabolic maximum principle gives
\begin{equation}\label{eq:u-bound}
 \min_M\cQ_F(\varphi_0)
 \leq u(x,t)\leq
 \max_M\cQ_F(\varphi_0),
\qquad 0\leq t<T_{\max}.
\end{equation}
Here we first apply the maximum principle on each compact time interval
$[0,T]$, $T<T_{\max}$, and then let $T\uparrow T_{\max}$.

Combining \eqref{eq:Q}, \eqref{eq:G-trace}, and
\eqref{eq:u-def}, we obtain
\begin{equation}\label{eq:inverse-trace-identity}
 \tr_{\chi_\varphi}\omega
 =n\exp\{F-u\}.
\end{equation}
Since $F$ is bounded on $M$ and $u$ satisfies \eqref{eq:u-bound}, there is a
constant $C\geq1$, depending only on the fixed data and the initial potential,
such that
\begin{equation}\label{eq:inverse-trace-bound}
 C^{-1}\leq\tr_{\chi_\varphi}\omega\leq C
\end{equation}
throughout the maximal interval.

\section{The scalar comparison principle}\label{sec:comparison}

We next establish the scalar order-preserving property underlying the barrier
argument.  The proof is a direct maximum-principle argument based only on the
monotonicity of the inverse-trace operator.  In particular, it is valid on a
general Hermitian background and does not require any closedness assumption on
$\omega$ or $\chi$.

\begin{proposition}[Parabolic comparison]\label{prop:comparison}
Let $\underline\varphi,\overline\varphi$ be smooth admissible functions on
$M\times[0,T]$ satisfying
\begin{align}
 \partial_t\underline\varphi
 &\leq \cQ_F(\underline\varphi),\label{eq:sub}\\
 \partial_t\overline\varphi
 &\geq \cQ_F(\overline\varphi).\label{eq:super}
\end{align}
If $\underline\varphi(\cdot,0)\leq\overline\varphi(\cdot,0)$, then
\[
 \underline\varphi\leq\overline\varphi
 \quad\text{on }M\times[0,T].
\]
\end{proposition}

\begin{proof}
For $\varepsilon>0$, define
\[
 w=\underline\varphi-\overline\varphi-\varepsilon t.
\]
If the conclusion fails, then $w$ has a positive maximum at some
$(x_0,t_0)\in M\times(0,T]$, since $w(\cdot,0)\leq0$.
At this point,
\[
 \partial_tw\geq0,
 \qquad
 \sqrt{-1}\,\partial\bar\partial w\leq0.
\]

Since $-\varepsilon t$ has no spatial derivatives,
\[
 \chi_{\underline\varphi}(x_0,t_0)
 \leq
 \chi_{\overline\varphi}(x_0,t_0).
\]
By the monotonicity of $G$ established in \eqref{eq:DG-expanded}, and since
the same function $F$ occurs in both operators,
\[
 \cQ_F(\underline\varphi)-\cQ_F(\overline\varphi)
 =
 G(\chi_{\underline\varphi})-G(\chi_{\overline\varphi})
 \leq0
\]
at $(x_0,t_0)$.

On the other hand, \eqref{eq:sub} and \eqref{eq:super} give
\[
 \varepsilon
 \leq\partial_t\underline\varphi-
 \partial_t\overline\varphi
 \leq\cQ_F(\underline\varphi)-
 \cQ_F(\overline\varphi)\leq0,
\]
a contradiction.  Letting $\varepsilon\downarrow0$ proves the result.
\end{proof}

\begin{corollary}[Oscillation nonexpansiveness]\label{cor:contraction}
Let $\varphi$ and $\psi$ be two solutions of \eqref{eq:flow} on $[0,T]$
with initial data $\varphi_0$ and $\psi_0$, respectively.  Then
\begin{align}
 \sup_M(\varphi(t)-\psi(t))
 &\leq\sup_M(\varphi(0)-\psi(0)),\label{eq:sup-contract}\\
 \inf_M(\varphi(t)-\psi(t))
 &\geq\inf_M(\varphi(0)-\psi(0)).\label{eq:inf-contract}
\end{align}
Consequently,
\begin{equation}\label{eq:osc-contract}
 \osc_M(\varphi(t)-\psi(t))
 \leq\osc_M(\varphi(0)-\psi(0)).
\end{equation}
\end{corollary}

\begin{proof}
The operator $\cQ_F$ is invariant under addition of spatial constants.
Therefore
\[
 \psi+\sup_M(\varphi_0-\psi_0)
\]
is a solution whose initial value lies above $\varphi_0$, while
\[
 \psi+\inf_M(\varphi_0-\psi_0)
\]
is a solution whose initial value lies below $\varphi_0$.
Applying Proposition \ref{prop:comparison} to the two pairs gives
\eqref{eq:sup-contract} and \eqref{eq:inf-contract}, and hence
\eqref{eq:osc-contract}.
\end{proof}

On the quotient $\cH_\chi/\mathbb R$, define
\begin{equation}\label{eq:osc-distance}
 \dosc([\varphi],[\psi])
 =
 \osc_M(\varphi-\psi).
\end{equation}
This is a well-defined metric.  Indeed, adding constants to either
representative does not change the oscillation, and
\[
 \osc_M(f+g)\leq\osc_M f+\osc_M g
\]
gives the triangle inequality.  Moreover,
$\dosc([\varphi],[\psi])=0$ if and only if $\varphi-\psi$ is constant.  Corollary \ref{cor:contraction} can therefore be restated as follows.

\begin{corollary}[Nonexpansive solution map]\label{cor:semiflow}
For $t\geq0$, let
\[
 \mathcal D_t
 :=
 \left\{
 [\varphi_0]\in\cH_\chi/\mathbb R:
 \text{the solution with initial value $\varphi_0$ exists on $[0,t]$}
 \right\}.
\]
The solution map
\[
 S_t:\mathcal D_t\longrightarrow\cH_\chi/\mathbb R,
 \qquad
 S_t[\varphi_0]=[\varphi(t)],
\]
is well defined and nonexpansive:
\begin{equation}\label{eq:semiflow-nonexpansive}
 \dosc(S_t[\varphi_0],S_t[\psi_0])
 \leq
 \dosc([\varphi_0],[\psi_0])
\end{equation}
for all $[\varphi_0],[\psi_0]\in\mathcal D_t$.
\end{corollary}

\begin{proof}
Since $\cQ_F(\varphi+c)=\cQ_F(\varphi)$ for every spatial constant $c$,
adding a constant to the initial potential adds the same constant to the
corresponding solution.  Hence both $\mathcal D_t$ and $S_t$ are well defined
on the quotient.  The estimate follows directly from Corollary
\ref{cor:contraction}.
\end{proof}

\begin{remark}\label{rem:comparison-hym}
A tensorial comparison principle for prescribed Hermitian--Yang--Mills flows
was recently developed in \cite{xyy26a}, where it is used to obtain uniform
$C^0$ estimates for Hermitian bundle metrics.  The scalar comparison argument
above is self-contained and uses only the monotonicity of $G$.  The structural
parallel is that ordered sub- and supersolutions provide pointwise barriers
for the evolving object, without requiring an energy functional.
\end{remark}

\section{Common-drift barriers and the cone condition}\label{sec:cone}

The two barrier inequalities have complementary roles.  Their combination
provides a time-independent $C^0$ estimate after subtraction of a linear drift,
while the lower barrier alone supplies the strict cone condition needed in the
second-order estimate.

\begin{lemma}[The lower barrier gives the cone condition]\label{lem:barrier-cone}
If $\underline\varphi\in\cH_\chi$ and
$\cQ_F(\underline\varphi)\geq\lambda$, then
\begin{equation}\label{eq:cone-from-barrier}
 n\chi_{\underline\varphi}^{n-1}
 >(n-1)e^{\lambda-F}
 \chi_{\underline\varphi}^{n-2}\wedge\omega.
\end{equation}
In particular, $[\chi]\in\cC_1(F-\lambda)$.
\end{lemma}

\begin{proof}
Fix $x\in M$.  After a complex linear change of holomorphic coordinates at
$x$, we may assume that
\[
 \chi_{\underline\varphi}(x)
 =
 \sqrt{-1}\sum_{i=1}^n dz^i\wedge d\bar z^i,
 \qquad
 \omega(x)
 =
 \sqrt{-1}\sum_{i=1}^n\mu_i\,dz^i\wedge d\bar z^i,
\]
where $\mu_i>0$.  Evaluating
$\cQ_F(\underline\varphi)\geq\lambda$ at $x$ gives
\begin{equation}\label{eq:trace-sub}
 \sum_{i=1}^n\mu_i
 =
 \tr_{\chi_{\underline\varphi}}\omega
 \leq
 ne^{F(x)-\lambda}.
\end{equation}

For each $i$, let $\Theta_i$ denote the standard positive
$(n-1,n-1)$-form obtained by omitting the $i$-th factor.  The coefficient of
$\Theta_i$ in
\[
 n\chi_{\underline\varphi}^{n-1}
 -(n-1)e^{\lambda-F}
 \chi_{\underline\varphi}^{n-2}\wedge\omega
\]
at $x$ is
\[
 (n-1)!
 \left(
 n-e^{\lambda-F(x)}
 \sum_{j\neq i}\mu_j
 \right).
\]
By \eqref{eq:trace-sub} and $\mu_i>0$,
\[
 e^{\lambda-F(x)}
 \sum_{j\neq i}\mu_j
 <
 e^{\lambda-F(x)}
 \sum_{j=1}^n\mu_j
 \leq n.
\]
Hence every coefficient is strictly positive.  Since $x$ was arbitrary,
\eqref{eq:cone-from-barrier} follows.
\end{proof}

\begin{lemma}[Common-drift trapping]\label{lem:trapping}
Assume \eqref{eq:barriers}.  For every $\varphi_0\in\cH_\chi$, there are
constants $A_0,B_0\in\mathbb R$ such that the solution of \eqref{eq:flow}
satisfies
\begin{equation}\label{eq:trapping}
 \underline\varphi+A_0
 \leq\varphi(t)-\lambda t
 \leq\overline\varphi+B_0
\end{equation}
throughout its maximal existence interval.  In particular,
\begin{equation}\label{eq:barrier-C0}
 \sup_{M\times[0,T_{\max})}
 |\varphi(x,t)-\lambda t|\leq C.
\end{equation}
\end{lemma}

\begin{proof}
Set
\[
 A_0:=\min_M(\varphi_0-\underline\varphi),
 \qquad
 B_0:=\max_M(\varphi_0-\overline\varphi).
\]
Then
\[
 \underline\varphi+A_0
 \leq\varphi_0
 \leq\overline\varphi+B_0
 \quad\text{on }M.
\]

The function
\[
 \underline\Phi(x,t)=\underline\varphi(x)+A_0+\lambda t
\]
satisfies
\[
 \partial_t\underline\Phi=\lambda
 \leq\cQ_F(\underline\varphi)
 =\cQ_F(\underline\Phi),
\]
whereas
\[
 \overline\Phi(x,t)=\overline\varphi(x)+B_0+\lambda t
\]
satisfies
\[
 \partial_t\overline\Phi=\lambda
 \geq\cQ_F(\overline\varphi)
 =\cQ_F(\overline\Phi).
\]

Since $\cQ_F$ is invariant under addition of functions that are constant in
space,
\[
 \cQ_F(\underline\Phi)=\cQ_F(\underline\varphi),
 \qquad
 \cQ_F(\overline\Phi)=\cQ_F(\overline\varphi).
\]
Thus $\underline\Phi$ and $\overline\Phi$ are respectively a subsolution and
a supersolution of \eqref{eq:flow}.

Fix any $T<T_{\max}$.  Proposition \ref{prop:comparison}, applied on
$M\times[0,T]$, gives
\[
 \underline\Phi\leq\varphi\leq\overline\Phi.
\]
Since $T<T_{\max}$ is arbitrary, \eqref{eq:trapping} holds throughout the
maximal existence interval.

Moreover,
\[
 |\varphi(x,t)-\lambda t|
 \leq
 \max\left\{
 \|\underline\varphi+A_0\|_{C^0(M)},
 \|\overline\varphi+B_0\|_{C^0(M)}
 \right\}
 =:C,
\]
which proves \eqref{eq:barrier-C0}.  In particular, $C$ is independent of
$t$ and of any terminal time $T<T_{\max}$.
\end{proof}

\begin{lemma}[Control of the quantity in Sun's estimate]\label{lem:sun-quantity}
Assume the hypotheses of Lemma~\ref{lem:trapping}, set
\[
 \psi=\varphi-\lambda t,\qquad
 \chi_0=\chi_{\varphi_0},\qquad
 v=\psi-\varphi_0.
\]
Then
\[
 v(\cdot,0)=0,
 \qquad
 \chi_0+\ddbar v=\chi_\varphi.
\]
Moreover, there is a constant $C_1>0$, independent of $t$ and of any
terminal time $T<T_{\max}$, such that
\begin{equation}\label{eq:sun-quantity-bound}
 0
 \leq
 v(x,t)-\inf_{M\times[0,t]}v
 \leq2C_1,
 \qquad
 (x,t)\in M\times[0,T_{\max}).
\end{equation}
\end{lemma}

\begin{proof}
By \eqref{eq:barrier-C0},
\[
 \|\psi(t)\|_{C^0(M)}\leq C
 \qquad\text{for all }t<T_{\max}.
\]
Hence
\[
 |v(x,t)|
 =
 |\psi(x,t)-\varphi_0(x)|
 \leq
 C+\|\varphi_0\|_{C^0(M)}
 =:C_1.
\]
For every $t<T_{\max}$,
\[
 -C_1
 \leq
 \inf_{M\times[0,t]}v
 \leq
 v(x,t)
 \leq C_1,
\]
which gives
\eqref{eq:sun-quantity-bound}.
\end{proof}

\begin{remark}[Why no space--time Moser iteration is needed]\label{rem:no-spacetime-moser}
The $C^0$ estimate used in this paper is \eqref{eq:barrier-C0}, obtained
entirely from the comparison principle.  No Moser iteration on
$M\times[0,T]$ is performed.  The argument of Tosatti--Weinkove
\cite{twasian} is elliptic and integrates on $M$ for a single function.  Applying it separately at each fixed time does not, without an additional
time-uniform estimate, control the space--time quantity
\[
 v(x,t)-\inf_{M\times[0,t]}v
\]
that appears in Sun's parabolic estimate.  The temporal information needed
there is supplied instead by the two-sided barrier estimate, as shown in
Lemma~\ref{lem:sun-quantity}.  Thus the parabolic zero-order argument is
logically independent of an extension of the Tosatti--Weinkove integration
argument to $M\times[0,T]$.
\end{remark}

\section{Second-order and higher-order estimates}\label{sec:estimates}

We now give the analytic continuation argument in detail.  The difficult local
second-order calculation, including the torsion terms on a Hermitian manifold,
is contained in Sun's estimate \cite{sun15}. We record the precise specialization needed here and then combine it with
the time-uniform drift bound obtained in Section~\ref{sec:cone}.

\begin{proposition}[Specialized parabolic second-order estimate]\label{prop:sun-c2}
Let $\widehat F\in C^\infty(M)$ and suppose that
$[\chi]\in\cC_1(\widehat F)$.  Let $v$ be an admissible solution of
\begin{equation}\label{eq:general-shifted}
 \partial_tv=
 \log\frac{(\chi_0+\ddbar v)^n}
 {\omega\wedge(\chi_0+\ddbar v)^{n-1}}
 +\widehat F,
 \qquad v(\cdot,0)=0,
\end{equation}
where $\chi_0>0$ and $[\chi_0]=[\chi]$.  Then there are constants
$A,C>0$, depending only on the fixed geometric data, $\widehat F$, the initial
form $\chi_0$, and a strict representative of the cone condition, such that
\begin{equation}\label{eq:sun-c2}
 \tr_\omega(\chi_0+\ddbar v)(x,t)
 \leq C\exp\left(A\left[v(x,t)-
 \inf_{M\times[0,t]}v\right]\right)
\end{equation}
for every $t$ in the existence interval.  The constants are independent of the
terminal time.
\end{proposition}

\begin{proof}
Set
\[
 \Psi=e^{-\widehat F}.
\]
Then \eqref{eq:general-shifted} takes the form
\[
 \partial_t v
 =
 \log\frac{(\chi_0+\ddbar v)^n}
 {(\chi_0+\ddbar v)^{n-1}\wedge\omega}
 -\log\Psi,
\]
while the cone condition becomes
\[
 n\widehat\chi^{n-1}
 >
 (n-1)\Psi\,
 \widehat\chi^{n-2}\wedge\omega.
\]
Thus this is the case $\alpha=1$ of Sun's parabolic complex
Monge--Amp\`ere type equation \cite[Theorem~1.1 and
Proposition~3.1]{sun15}.

We note that \cite{sun15} uses the convention
$\chi_u=\chi+\frac{\sqrt{-1}}{2}\partial\bar\partial u$.
The change of variables
\[
 u(x,s)=2v(x,s/2)
\]
matches that convention with ours and leaves the parabolic equation in the
same form.  Sun's estimate therefore gives, after changing the constants
$A$ and $C$ if necessary,
\[
 \tr_\omega(\chi_0+\ddbar v)(x,t)
 \leq
 C\exp\left(
 A\left[
 v(x,t)-\inf_{M\times[0,t]}v
 \right]\right),
\]
which is \eqref{eq:sun-c2}.  The constants are independent of the terminal
time.
\end{proof}

\begin{proposition}[Uniform ellipticity from a drift bound]\label{prop:ellipticity}
Let $\varphi$ solve \eqref{eq:flow} on $[0,T)$ and suppose that for some
$\lambda\in\mathbb R$,
\begin{equation}\label{eq:drift-C0-assumption}
 \sup_{M\times[0,T)}|\varphi(x,t)-\lambda t|\leq C_0,
\end{equation}
and $[\chi]\in\cC_1(F-\lambda)$.  Then there is a constant $C\geq1$,
independent of $T$, such that
\begin{equation}\label{eq:metric-equivalence}
 C^{-1}\omega\leq\chi_\varphi(t)\leq C\omega
\end{equation}
for all $t<T$.
\end{proposition}

\begin{proof}
Set
\[
 \psi=\varphi-\lambda t.
\]
Since $\chi_\psi=\chi_\varphi$, the shifted function satisfies
\begin{equation}\label{eq:shifted-flow-detailed}
 \partial_t\psi=
 \log\frac{\chi_\psi^n}
 {\omega\wedge\chi_\psi^{n-1}}+F-\lambda.
\end{equation}
Write
\[
 \chi_0=\chi_{\varphi_0},
 \qquad
 v=\psi-\varphi_0.
\]
Then $v(\cdot,0)=0$, $\chi_0+\ddbar v=\chi_\psi$, and $v$ solves
\eqref{eq:general-shifted} with $\widehat F=F-\lambda$.
Proposition \ref{prop:sun-c2} gives
\begin{equation}\label{eq:upper-trace-calculation}
 \tr_\omega\chi_\varphi(x,t)
 \leq C\exp\left(A\left[v(x,t)-
 \inf_{M\times[0,t]}v\right]\right).
\end{equation}

By the drift bound \eqref{eq:drift-C0-assumption} and the fixed initial
potential $\varphi_0$, the argument of Lemma~\ref{lem:sun-quantity} gives
\[
 0
 \leq
 v(x,t)-\inf_{M\times[0,t]}v
 \leq 2C_1,
\]
where
\[
 C_1=C_0+\|\varphi_0\|_{C^0(M)}
\]
is independent of $t$ and $T$.  Hence
\eqref{eq:upper-trace-calculation} yields
\begin{equation}\label{eq:trace-upper}
 \tr_\omega\chi_\varphi
 \leq
 Ce^{2AC_1}
 =:C_2,
\end{equation}
where $C_2$ is independent of the terminal time $T$.

On the other hand, \eqref{eq:inverse-trace-identity} and
\eqref{eq:u-bound} give
\begin{equation}\label{eq:inverse-trace-upper}
 \tr_{\chi_\varphi}\omega\leq C.
\end{equation}
Let $\mu_1,\ldots,\mu_n>0$ be the eigenvalues of $\chi_\varphi$ with
respect to $\omega$.  Then \eqref{eq:trace-upper} and
\eqref{eq:inverse-trace-upper} give
\[
 \sum_{i=1}^n\mu_i\leq C,
 \qquad
 \sum_{i=1}^n\mu_i^{-1}\leq C.
\]
In particular,
\[
 C^{-1}\leq\mu_i\leq C,
 \qquad 1\leq i\leq n,
\]
after enlarging $C$ if necessary.  This proves
\eqref{eq:metric-equivalence}.
\end{proof}

\begin{proposition}[Higher-order estimates]\label{prop:higher}
Under the hypotheses of Proposition \ref{prop:ellipticity}, for every
$\tau>0$ and every integer $m\geq0$, there exists a constant
$C_{m,\tau}$, depending only on $\tau$, the fixed data, the initial
potential, and the constants in Proposition \ref{prop:ellipticity}, but not
on the terminal time $T$, such that whenever $T>\tau$,
\begin{equation}\label{eq:higher}
 \sup_{\tau\leq t<T}
 \|\varphi(t)-\lambda t\|_{C^m(M)}
 \leq C_{m,\tau}.
\end{equation}
In particular, the solution extends smoothly past every finite terminal time.
\end{proposition}

\begin{proof}
Let $\psi=\varphi-\lambda t$.  By Proposition
\ref{prop:ellipticity}, the eigenvalues of $\chi_\psi$ with respect to $\omega$
remain in a fixed compact subset of the positive cone.  Formula
\eqref{eq:DG-expanded} therefore gives uniform ellipticity:
there are constants $0<\theta\leq\Theta$ such that, in every fixed coordinate
chart,
\begin{equation}\label{eq:uniform-parabolicity}
 \theta |\xi|^2\leq
 a^{i\bar j}(x,t)\xi_i\overline{\xi_j}
 \leq\Theta |\xi|^2.
\end{equation}

The equation is therefore uniformly parabolic and concave in the complex
Hessian.  Standard parabolic Evans--Krylov regularity
\cite{evans82,krylov84}, in the complex-Hessian setting used in
\cite{sun15,pt21}, gives
some $\alpha\in(0,1)$ such that, for every $\tau>0$ with $T>\tau$,
\begin{equation}\label{eq:EK}
 \sup_{\tau\leq t<T}
 \left(
 \|\psi(t)\|_{C^{2,\alpha}(M)}
 +
 \|\partial_t\psi(t)\|_{C^\alpha(M)}
 \right)
 \leq C_\tau.
\end{equation}
Here $C_\tau$ is independent of $T$.  This is the standard regularity step
used after the uniform $C^0$ and $C^2$ estimates in
\cite{sun15,pt21}.  Equivalently, one may apply the local parabolic
Evans--Krylov estimate on backward cylinders of fixed size contained in
$M\times(0,T)$.

We stress that this is a complex-Hessian regularity statement; no uniform
ellipticity in the complementary directions of the full real Hessian is
being asserted.

We briefly indicate the subsequent Schauder bootstrap.  Differentiating
\eqref{eq:shifted-flow-detailed} in a spatial direction $D$ yields
\begin{equation}\label{eq:first-diff}
 (\partial_t-L_t)D\psi
 =D(F-\lambda)+\mathcal E_1,
\end{equation}
where $\mathcal E_1$ contains only derivatives of the fixed background forms and
coefficients controlled by \eqref{eq:EK}.  Linear parabolic Schauder estimates
give a $C^{3+\alpha}$ bound.  Repeated differentiation produces equations of the
form
\begin{equation}\label{eq:bootstrap-eq}
 (\partial_t-L_t)D^r\psi=\mathcal E_r,
\end{equation}
where $\mathcal E_r$ depends on derivatives of $\psi$ of order at most $r+1$
and on fixed smooth data.  Induction and Schauder theory give
\eqref{eq:higher} for every $m$.

It remains to justify the continuation statement.  Suppose that $T<\infty$.
Applying the preceding estimates with $\tau=T/2$, we obtain uniform spatial
$C^m$ bounds for $\psi$ on $[T/2,T)$ for every $m$.  Since
\[
 \partial_t\psi
 =
 \cQ_{F-\lambda}(\psi)
\]
is a smooth function of the fixed data and of the complex Hessian of
$\psi$, these bounds also imply, for every $m\geq0$,
\[
 \sup_{T/2\leq t<T}
 \|\partial_t\psi(t)\|_{C^m(M)}
 \leq C_m.
\]
Consequently, for $T/2\leq t_1<t_2<T$,
\[
 \|\psi(t_2)-\psi(t_1)\|_{C^m(M)}
 \leq
 \int_{t_1}^{t_2}
 \|\partial_t\psi(s)\|_{C^m(M)}\,ds
 \leq
 C_m|t_2-t_1|.
\]
Thus $\psi(t)$ is Cauchy in every $C^m(M)$ as $t\nearrow T$ and converges
smoothly to a function $\psi_T$.  The metric equivalence
\eqref{eq:metric-equivalence} passes to the limit, so
\[
 \chi_{\psi_T}>0.
\]
Hence $\psi_T$ is admissible.  Thus the solution extends smoothly past $T$.
In particular, if $T=T_{\max}<\infty$, this contradicts the maximality of
$T_{\max}$.  Therefore no finite maximal existence time can occur while the
hypotheses of Proposition \ref{prop:ellipticity} hold uniformly.
\end{proof}

\section{A uniform Harnack inequality and exponential convergence}\label{sec:convergence}

We next derive exponential convergence from a standard parabolic Harnack
inequality.  Once uniform ellipticity and higher-order bounds are available,
the Harnack constant is uniform on all time translates.

\begin{lemma}[Uniform Harnack inequality]\label{lem:harnack}
Let $w\geq0$ be a smooth solution of
\begin{equation}\label{eq:linear-w}
 \partial_t w
 =
 a^{i\bar j}(x,t)w_{i\bar j}
\end{equation}
on $M\times[s,s+1]$.  Suppose that the coefficient matrix
$(a^{i\bar j})$ satisfies the uniform ellipticity bounds
\eqref{eq:uniform-parabolicity}.  Then there exists a constant
$C_H>1$, depending only on the uniform ellipticity constants and the
fixed background geometry, but not on $s$, such that
\begin{equation}\label{eq:harnack}
 \sup_M w(\cdot,s+\tfrac12)
 \leq
 C_H\inf_M w(\cdot,s+1).
\end{equation}
\end{lemma}

\begin{proof}
In local real coordinates, the operator
$a^{i\bar j}\partial_i\partial_{\bar j}$ is a uniformly elliptic
nondivergence-form operator, with ellipticity constants controlled by
\eqref{eq:uniform-parabolicity}.  The parabolic Krylov--Safonov Harnack
inequality therefore applies in each fixed coordinate cylinder; see, for
example, \cite[Chapter~VII]{lieb}.  A finite covering of the compact connected
manifold $M$, together with a finite Harnack chain, yields
\eqref{eq:harnack}.  All geometric and ellipticity constants are uniform, so
$C_H$ is independent of the time translate $s$.  The estimate for
nonnegative solutions follows by applying the positive-solution estimate to
$w+\varepsilon$ and letting $\varepsilon\downarrow0$.
\end{proof}

\begin{proposition}[Exponential decay of the speed oscillation]\label{prop:speed-osc}
Assume that the solution exists on $[0,\infty)$ and that the estimates
\eqref{eq:metric-equivalence} and \eqref{eq:higher} hold uniformly on
$[1,\infty)$.  Let $u=\partial_t\varphi$.  Then there are constants
$C,\eta>0$ such that
\begin{equation}\label{eq:osc-decay}
 \osc_Mu(\cdot,t)\leq Ce^{-\eta t}
\end{equation}
for all $t\geq1$.
\end{proposition}

\begin{proof}
Set
\[
 M(t)=\sup_Mu(\cdot,t),
 \qquad
 m(t)=\inf_Mu(\cdot,t),
 \qquad
 \Theta(t)=M(t)-m(t).
\]
By \eqref{eq:u-evolution} and the maximum principle, $M$ is nonincreasing and
$m$ is nondecreasing.  Fix an integer $j\geq1$.  On $M\times[j,j+1]$ the
functions
\[
 v_j(x,t)=M(j)-u(x,t),
 \qquad
 w_j(x,t)=u(x,t)-m(j).
\]
Since $M(t)$ is nonincreasing and $m(t)$ is nondecreasing,
\[
 u(x,t)\leq M(t)\leq M(j),
 \qquad
 u(x,t)\geq m(t)\geq m(j)
\]
for $t\in[j,j+1]$.  Hence $v_j,w_j\geq0$ on
$M\times[j,j+1]$.  Moreover, since constants are annihilated by $L_t$,
both $v_j$ and $w_j$ solve \eqref{eq:linear-w}.

 Applying Lemma
\ref{lem:harnack} to $v_j+\varepsilon$ and then letting
$\varepsilon\downarrow0$ yields
\begin{align}
 M(j)-m(j+\tfrac12)
 &=\sup_Mv_j(\cdot,j+\tfrac12)\notag\\
 &\leq C_H\inf_Mv_j(\cdot,j+1)\notag\\
 &=C_H\bigl(M(j)-M(j+1)\bigr).
 \label{eq:harnack-v}
\end{align}
The same argument for $w_j$ gives
\begin{align}
 M(j+\tfrac12)-m(j)
 &\leq C_H\bigl(m(j+1)-m(j)\bigr).
 \label{eq:harnack-w}
\end{align}
Adding \eqref{eq:harnack-v} and \eqref{eq:harnack-w}, we obtain
\begin{equation}\label{eq:theta-discrete-1}
 \Theta(j)+\Theta(j+\tfrac12)
 \leq C_H\bigl(\Theta(j)-\Theta(j+1)\bigr).
\end{equation}
Dropping the nonnegative middle-time term gives
\[
 \Theta(j+1)\leq q\Theta(j),
 \qquad
 q:=1-C_H^{-1}\in(0,1).
\]
Thus
\[
 \Theta(j)\leq q^{j-1}\Theta(1),
 \qquad j\geq1.
\]
Set
\[
 \eta:=-\log q>0.
\]
For $t\geq1$, let $j=\lfloor t\rfloor$.  Since $\Theta$ is nonincreasing,
\[
 \Theta(t)
 \leq
 \Theta(j)
 \leq
 q^{j-1}\Theta(1)
 \leq
 C e^{-\eta t}.
\]
This proves \eqref{eq:osc-decay}.
\end{proof}

\begin{proposition}[Convergence under a uniform drift bound]\label{prop:drift-convergence}
Let $\varphi$ be a solution on its maximal interval $[0,T_{\max})$.
Suppose that, for some $\lambda\in\mathbb R$,
\begin{equation}\label{eq:drift-c0}
 \sup_{M\times[0,T_{\max})}|\varphi(x,t)-\lambda t|<\infty
\end{equation}
and $[\chi]\in\cC_1(F-\lambda)$.  Then $T_{\max}=\infty$.  Moreover, there
is a smooth $\varphi_\infty\in\cH_\chi$ satisfying
\begin{equation}\label{eq:limit-equation}
 \cQ_F(\varphi_\infty)=\lambda
\end{equation}
such that, for every integer $k\geq0$, there are constants
$C_k,\eta_k>0$ with
\begin{equation}\label{eq:drift-exp}
 \|\varphi(t)-\lambda t-\varphi_\infty\|_{C^k(M)}
 \leq C_ke^{-\eta_kt}.
\end{equation}
\end{proposition}

\begin{proof}
Propositions \ref{prop:ellipticity} and \ref{prop:higher} give uniform
ellipticity, all higher-order estimates, and $T_{\max}=\infty$.
Proposition \ref{prop:speed-osc} gives
\begin{equation}\label{eq:u-osc-exp}
 M(t)-m(t)\leq Ce^{-\eta t}.
\end{equation}

Since $M$ is nonincreasing and $m$ is nondecreasing, both have limits.
By \eqref{eq:u-osc-exp}, these limits coincide; denote their common value by
$c$.  Hence
\[
 m(t)\leq c\leq M(t),
\]
and therefore
\[
 \|u(\cdot,t)-c\|_{C^0(M)}
 \leq
 M(t)-m(t)
 \leq
 Ce^{-\eta t}.
\]

Set
\[
 \psi=\varphi-\lambda t,
 \qquad
 V=\int_M\omega^n,
 \qquad
 \overline\psi(t)=\frac{1}{V}\int_M\psi(t)\,\omega^n.
\]

By \eqref{eq:drift-c0}, $\overline\psi(t)$ is uniformly bounded.  Moreover,
\begin{equation}\label{eq:mean-derivative}
 \overline\psi'(t)
 =
 \frac{1}{V}\int_M(u-\lambda)\,\omega^n
 \longrightarrow
 c-\lambda.
\end{equation}
Suppose that $c\neq\lambda$.  Then there exists $t_0>0$ such that
\[
 |\overline\psi'(t)|
 \geq
 \frac{|c-\lambda|}{2},
 \qquad t\geq t_0,
\]
and $\overline\psi'(t)$ has the sign of $c-\lambda$.  Consequently,
\[
 |\overline\psi(t)-\overline\psi(t_0)|
 \geq
 \frac{|c-\lambda|}{2}(t-t_0),
 \qquad t\geq t_0,
\]
contradicting the uniform boundedness of $\overline\psi$.  Thus
$c=\lambda$.

Since $m(t)\leq c=\lambda\leq M(t)$, \eqref{eq:u-osc-exp} now gives
\begin{equation}\label{eq:psi-t-c0}
 \|\partial_t\psi(t)\|_{C^0(M)}
 =
 \|u(t)-\lambda\|_{C^0(M)}
 \leq
 Ce^{-\eta t}.
\end{equation}

To upgrade the decay, fix $k\geq0$ and choose an integer $r>k$.
By Proposition \ref{prop:higher}, applied with sufficiently high order,
the quantities
\[
 \|\psi(t)\|_{C^{r+2}(M)}
\]
are uniformly bounded for $t\geq1$.  Since
\[
 \partial_t\psi
 =
 \cQ_{F-\lambda}(\psi)
\]
depends smoothly on the fixed data and on the complex Hessian of $\psi$,
while $\chi_\psi$ remains uniformly equivalent to $\omega$, it follows that
\[
 \sup_{t\geq1}
 \|\partial_t\psi(t)\|_{C^r(M)}
 \leq C_r.
\]
The interpolation inequality therefore gives, for some
$\theta=\theta(k,r)\in(0,1)$,
\begin{align}
 \|\partial_t\psi(t)\|_{C^k}
 &\leq
 C\|\partial_t\psi(t)\|_{C^0}^{1-\theta}
  \|\partial_t\psi(t)\|_{C^r}^{\theta}\notag\\
 &\leq
 C_ke^{-\eta_k t}.
 \label{eq:psi-t-ck}
\end{align}
Hence, for $t_2>t_1\geq1$,
\[
 \|\psi(t_2)-\psi(t_1)\|_{C^k}
 \leq\int_{t_1}^{t_2}
 \|\partial_s\psi(s)\|_{C^k}\,ds
 \leq C_ke^{-\eta_k t_1}.
\]

Thus $\psi(t)$ is Cauchy in every $C^k$ norm and converges smoothly to a
function $\varphi_\infty$. The uniform metric equivalence \eqref{eq:metric-equivalence} passes to the limit, so $\chi_{\varphi_\infty}>0$ and hence $\varphi_\infty\in\cH_\chi$. Letting $t\to\infty$ in
\eqref{eq:shifted-flow-detailed} and using
$\partial_t\psi\to0$ proves \eqref{eq:limit-equation}; integrating
\eqref{eq:psi-t-ck} from $t$ to infinity gives \eqref{eq:drift-exp}.
\end{proof}

\begin{corollary}[Asymptotic synchronization]\label{cor:synchronization}
Assume that \eqref{eq:elliptic} is solvable.  Let $\varphi$ and $\psi$ be the
solutions of \eqref{eq:flow} starting from arbitrary elements of $\cH_\chi$.
Then there exists a constant $c_\infty\in\mathbb R$ such that, for every
integer $k\geq0$,
\begin{equation}\label{eq:synchronization}
 \|\varphi(t)-\psi(t)-c_\infty\|_{C^k(M)}
 \leq C_ke^{-\eta_kt}.
\end{equation}
In particular, the elliptic solution class is the unique globally attracting
fixed point of the induced semiflow on $\cH_\chi/\mathbb R$.
\end{corollary}

\begin{proof}
By Theorem \ref{thm:barrier}, solvability of \eqref{eq:elliptic} implies
global existence and exponential convergence for every admissible initial
potential. Let $\varphi_\infty$ and $\psi_\infty$ be the corresponding
drift-corrected limits. Proposition \ref{prop:elliptic-unique} gives
\[
\varphi_\infty-\psi_\infty=c_\infty
\]
for some constant $c_\infty$. Subtracting the two exponential convergence
estimates yields \eqref{eq:synchronization}.

Since global existence holds for
every initial class, the solution maps define a global semiflow on
$\cH_\chi/\mathbb R$; its nonexpansiveness follows from Corollary
\ref{cor:semiflow}, and the uniqueness of the elliptic solution class gives
the final assertion.
\end{proof}

\begin{proposition}[Uniqueness of the elliptic pair]\label{prop:elliptic-unique}
If $(\varphi_1,\lambda_1)$ and $(\varphi_2,\lambda_2)$ solve
\eqref{eq:elliptic}, then
\(
 \lambda_1=\lambda_2,
\)
and $\varphi_1-\varphi_2$ is constant.
\end{proposition}

\begin{proof}
Let $x_{\max}$ be a maximum point of
$\varphi_1-\varphi_2$.  Then
\[
 \ddbar(\varphi_1-\varphi_2)(x_{\max})\leq0,
\]
and hence
\[
 \chi_{\varphi_1}(x_{\max})
 \leq
 \chi_{\varphi_2}(x_{\max}).
\]
Since $G$ is increasing and the same function $F$ occurs in both elliptic
equations,
\[
 \lambda_1
 =
 G(\chi_{\varphi_1})+F
 \leq
 G(\chi_{\varphi_2})+F
 =
 \lambda_2
\]
at $x_{\max}$.  Applying the same argument at a minimum point gives
$\lambda_1\geq\lambda_2$.  Thus
\[
 \lambda_1=\lambda_2.
\]

Set
\[
 w=\varphi_1-\varphi_2,
 \qquad
 A_s=\chi_{\varphi_2}+s\ddbar w,
 \qquad 0\leq s\leq1.
\]
Since
\[
 A_s=(1-s)\chi_{\varphi_2}+s\chi_{\varphi_1},
\]
the convexity of the positive cone gives $A_s>0$.  As the two elliptic
equations have the same constant $\lambda_1=\lambda_2$, the function $F$
cancels and
\[
 0
 =
 G(\chi_{\varphi_1})-G(\chi_{\varphi_2})
 =
 \int_0^1DG_{A_s}(\ddbar w)\,ds
 =
 b^{i\bar j}w_{i\bar j},
\]
where
\[
 b^{i\bar j}
 =
 \int_0^1a_{A_s}^{i\bar j}\,ds.
\]
The coefficients $b^{i\bar j}$ are smooth and uniformly positive definite:
indeed, $A_s>0$ on the compact set $M\times[0,1]$.  Hence the strong maximum
principle implies that $w$ is constant, since $M$ is connected.
\end{proof}

\section{Proofs of the barrier results}\label{sec:barrier-proof}

\begin{proof}[Proof of Theorem \ref{thm:barrier}]
Lemma \ref{lem:barrier-cone} gives
\[
 [\chi]\in\cC_1(F-\lambda),
\]
while Lemma \ref{lem:trapping} gives
\[
 \sup_{M\times[0,T_{\max})}
 |\varphi(x,t)-\lambda t|<\infty.
\]
Thus the two hypotheses of Proposition \ref{prop:drift-convergence} are
satisfied, and that proposition yields $T_{\max}=\infty$ together with
exponential smooth convergence to a solution of
\[
 \cQ_F(\varphi_\infty)=\lambda.
\]
Proposition \ref{prop:elliptic-unique} gives uniqueness of the elliptic
solution modulo constants, while uniqueness of the parabolic solution follows
from Proposition \ref{prop:comparison}.
\end{proof}

\begin{proof}[Proof of Corollary \ref{cor:stability}]
Let $(\varphi_*,\lambda)$ be a solution of \eqref{eq:elliptic}.  Taking
\[
 \underline\varphi=\overline\varphi=\varphi_*
\]
in Theorem \ref{thm:barrier} gives global existence and exponential
convergence from every initial potential in $\cH_\chi$.  The synchronization
statement follows from Corollary \ref{cor:synchronization}.
\end{proof}

\begin{proof}[Proof of Corollary \ref{cor:equivalence}]
For (1)$\Rightarrow$(2), take the same elliptic solution as both barriers.
Theorem \ref{thm:barrier} gives (2)$\Rightarrow$(3).

It remains to prove (3)$\Rightarrow$(1).  Fix any
$\varphi_0\in\cH_\chi$ and set
\[
 \psi(t)=\varphi(t)-\lambda t.
\]
By assumption,
\[
 \psi(t)\longrightarrow\psi_\infty
 \quad\text{in }C^\infty(M)
\]
for some $\psi_\infty\in\cH_\chi$.  Since
\[
 \partial_t\psi=\cQ_{F-\lambda}(\psi),
\]
the smooth convergence and the admissibility of $\psi_\infty$ imply
\[
 \partial_t\psi
 \longrightarrow
 \cQ_{F-\lambda}(\psi_\infty)
 \quad\text{uniformly on }M.
\]
The limit on the right must vanish.  Indeed, if it were nonzero at some
point, then by continuity $\partial_t\psi$ would eventually have a fixed sign
and be bounded away from zero at that point, contradicting the convergence of
$\psi(t)$.  Hence
\[
 \cQ_{F-\lambda}(\psi_\infty)=0,
\]
or equivalently
\[
 \cQ_F(\psi_\infty)=\lambda.
\]
Thus (1) holds.
\end{proof}

\begin{remark}\label{rem:roles}
The lower and upper inequalities in \eqref{eq:barriers} are not redundant.  The
lower barrier yields the strict cone condition through Lemma
\ref{lem:barrier-cone}; the pair of barriers yields the uniform drift-corrected
$C^0$ estimate.  These are precisely the two inputs used in Proposition
\ref{prop:drift-convergence}.
\end{remark}

\section{Affine renormalization and an existence criterion}\label{sec:affine-proof}

We finish with the affine-renormalization criterion.  This argument produces one
elliptic solution by applying Sun's convergence theorem to a specially chosen
initial potential; the barrier theorem then upgrades that one solution to global
stability from every admissible initial value.

\begin{lemma}[Affine normalization]\label{lem:affine}
Let $\varphi$ solve \eqref{eq:flow}, fix $\kappa\in\mathbb R$, and set
$\psi=\varphi-\kappa t$.  Then
\begin{equation}\label{eq:affine-shift}
 \partial_t\psi=\cQ_{F-\kappa}(\psi).
\end{equation}
If
\begin{equation}\label{eq:initial-speed-sign}
 \kappa\geq\max_M\cQ_F(\varphi_0),
\end{equation}
then $\partial_t\psi\leq0$ throughout the maximal existence interval.
\end{lemma}

\begin{proof}
Because $\chi_\psi=\chi_\varphi$,
\[
 \partial_t\psi
 =\partial_t\varphi-\kappa
 =\cQ_F(\varphi)-\kappa
 =\cQ_{F-\kappa}(\psi).
\]

The speed $q=\partial_t\psi$ satisfies
\[
 \partial_tq=L_tq,
 \qquad
 q(\cdot,0)=\cQ_F(\varphi_0)-\kappa\leq0.
\]
Fix $T<T_{\max}$.  Since the solution is smooth and admissible on
$M\times[0,T]$, the operator $L_t$ is uniformly parabolic there.  The
parabolic maximum principle gives
\[
 q\leq0
 \quad\text{on }M\times[0,T].
\]
Since $T<T_{\max}$ is arbitrary, $\partial_t\psi=q\leq0$ throughout the
maximal existence interval.
\end{proof}

\begin{proof}[Proof of Theorem \ref{thm:affine}]
Start the flow at $\varphi_*$.  After replacing the background form $\chi$ by
$\chi_{\varphi_*}$ and the unknown by
\[
 v=\varphi-\varphi_*-\kappa t,
\]
we obtain
\begin{equation}\label{eq:affine-sun-form}
 \partial_tv=
 \log\frac{(\chi_{\varphi_*}+\ddbar v)^n}
 {\omega\wedge(\chi_{\varphi_*}+\ddbar v)^{n-1}}
 +F-\kappa,
 \qquad v(\cdot,0)=0.
\end{equation}

At $t=0$, condition \eqref{eq:affine-speed} is equivalent to
\begin{equation}\label{eq:sun-sign-match}
 \frac{\chi_{\varphi_*}^n}
 {\omega\wedge\chi_{\varphi_*}^{n-1}}
 \leq e^{\kappa-F}.
\end{equation}
Moreover, Lemma \ref{lem:affine} gives
\(
 \partial_t v\leq0
\)
throughout the maximal existence interval.

Assumption \eqref{eq:affine-cone} is precisely the cone condition for
\eqref{eq:affine-sun-form}, while \eqref{eq:sun-sign-match} is the
pointwise initial condition in \cite[Theorem~1.2]{sun15}.  After the harmless
normalization of the $\sqrt{-1}\partial\bar\partial$ convention described in
the proof of Proposition \ref{prop:sun-c2}, Sun's theorem therefore applies.

Set
\[
 \widetilde v(t)
 =
 v(t)-\frac{\int_Mv(t)\,\omega^n}{\int_M\omega^n}.
\]
Then $\widetilde v(t)$ converges smoothly to a function
$v_\infty$, and Sun's theorem gives a real constant $b$ such that
\[
 (\chi_{\varphi_*}+\ddbar v_\infty)^n
 =
 e^b e^{\kappa-F}
 (\chi_{\varphi_*}+\ddbar v_\infty)^{n-1}\wedge\omega.
\]
Define
\[
 \varphi_\infty:=\varphi_*+v_\infty,
 \qquad
 \lambda:=\kappa+b.
\]
Then
\[
 \frac{\chi_{\varphi_\infty}^n}
 {\omega\wedge\chi_{\varphi_\infty}^{n-1}}
 =
 e^{\lambda-F},
\]
and hence
\[
 \cQ_F(\varphi_\infty)=\lambda.
\]
Proposition \ref{prop:elliptic-unique} shows that $\lambda$ is the unique
constant for which \eqref{eq:elliptic} is solvable.  Corollary
\ref{cor:stability} then yields global existence and exponential convergence
from every initial potential in $\cH_\chi$.
\end{proof}

\begin{proof}[Proof of Corollary \ref{cor:hermitian}]
Take $\kappa=0$ in Theorem \ref{thm:affine}.  Condition
\eqref{eq:hermitian-super} is equivalent to
$\max_M\cQ_F(\varphi_*)\leq0$, while
\eqref{eq:hermitian-cone} is exactly \eqref{eq:affine-cone} with
$\kappa=0$.  The conclusion follows from Theorem \ref{thm:affine} and
Proposition \ref{prop:elliptic-unique}.
\end{proof}

\begin{corollary}\label{cor:pointwise}
Suppose that $\varphi_*\in\cH_\chi$ and $\kappa\in\mathbb R$ satisfy
\eqref{eq:affine-speed} and
\begin{equation}\label{eq:pointwise-cone-restated}
 \chi_{\varphi_*}-\frac{n-1}{n}e^{\kappa-F}\omega>0.
\end{equation}
Then all conclusions of Theorem \ref{thm:affine} hold.
\end{corollary}

\begin{proof}
Since $\chi_{\varphi_*}>0$, wedging the positive $(1,1)$-form in
\eqref{eq:pointwise-cone-restated} with the positive form
$n\chi_{\varphi_*}^{n-2}$ gives
\[
 n\chi_{\varphi_*}^{n-1}
 >
 (n-1)e^{\kappa-F}
 \chi_{\varphi_*}^{n-2}\wedge\omega.
\]
Thus $\chi_{\varphi_*}$ is a strict representative of
$[\chi]\in\cC_1(F-\kappa)$.
\end{proof}

\begin{remark}\label{rem:scope}
Theorem \ref{thm:affine} is a sufficient solvability criterion obtained from
Sun's theorem; it is not claimed to be necessary.  In particular, neither
Theorem \ref{thm:affine} nor Corollary \ref{cor:hermitian} asserts that the cone
condition alone gives convergence of the unshifted flow on a general Hermitian
manifold.
The sign condition is imposed only after the affine change of variables, and
its role is to place the special initial flow within the scope of Sun's
convergence theorem.  Once one elliptic solution has been obtained, the common-
drift barrier theorem removes the sign restriction for all other initial data.
\end{remark}

\section*{Acknowledgements}

The author would like to thank Professor Kefeng Liu for his guidance and help,
and Professor Yi Li for his interest in this work. The author also used
ChatGPT (OpenAI) to assist with language refinement, organization of the
exposition, and consistency checking. All mathematical arguments and final
content were reviewed and verified by the author.

\section*{Statements and Declarations}

\subsection*{Funding}
No funding was received to assist with the preparation of this manuscript.

\subsection*{Data availability}
No datasets were generated or analysed during the current study.

\subsection*{Competing interests}
The author declares no competing interests.

\end{document}